\documentclass[titlepage,11pt]{article}
\usepackage{latexsym}
\usepackage{amsfonts}
\usepackage{amsmath}
\usepackage{mathtools}
\usepackage{comment}
\usepackage{tikz}
\usepackage{float}
\newcommand\blackslug{\hbox{\hskip 1pt \vrule width 4pt height 8pt depth 1.5pt
        \hskip 1pt}}
\newcommand\bbox{\quad \blackslug \bigbreak}
\def\DD{\hbox{-}}
\def\CC{\hbox{-}\cdots\hbox{-}}
\def\LL{,\ldots,}
\newcommand{\vare}{\varepsilon}
\newcommand{\vep}{\varepsilon}
\newcommand{\erh}{Erd\H{o}s-Hajnal}
\newcommand{\ellh}{h}

\newcommand{\mab}{\mathbb}

\DeclarePairedDelimiter\floor{\lfloor}{\rfloor}%

\DeclarePairedDelimiter\abs{\lvert}{\rvert}%

\title{Induced subgraph density. V. All paths approach Erd\H{o}s-Hajnal}
\author{
Tung Nguyen\thanks{Supported by AFOSR grants
A9550-19-1-0187 and FA9550-22-1-0234, and by NSF grants  DMS-1800053 and DMS-2154169.}\\
Princeton University,\\ Princeton, NJ 08544, USA
\and
Alex Scott\thanks{Supported by EPSRC grant EP/X013642/1}\\
University of Oxford, \\
Oxford OX2 6GG, UK
\and
Paul Seymour\thanks{Supported by AFOSR grants
A9550-19-1-0187 and FA9550-22-1-0234, and by NSF grants  DMS-1800053 and DMS-2154169.}\\
Princeton University,\\ Princeton, NJ 08544, USA}

\date{July 4, 2022; revised \today}

\newtheorem{thm}{}[section]

\newcommand{\Proof}{\noindent{\bf Proof.}\ \ }

\begin{document}
\maketitle
\begin{abstract}
The Erd\H{o}s-Hajnal conjecture says that, for every graph $H$, there exists $c>0$ such that every $H$-free 
graph on $n$ vertices 
has a clique or stable set of size at least $n^c$.  In this paper we are concerned with the case when $H$ is a path.
The conjecture has been proved for paths with at most five vertices, but not for longer paths.
We prove that the conjecture is ``nearly'' true for all paths: 
for every path $H$, all $H$-free graphs with $n$ vertices have cliques or stable sets of size 
at least $2^{(\log n)^{1-o(1)}}$.
\end{abstract}

\section{Introduction}

A graph is {\em $H$-free} if it has no induced subgraph isomorphic to $H$.
A famous conjecture of Erd\H{o}s and Hajnal~\cite{EH77, EH89} from 1977 says:
\begin{thm}\label{EHconj}
{\bf The \erh{} Conjecture:} For every graph $H$ there exists $c>0$ such that every $H$-free graph $G$ has a stable set or
clique of size at least $|G|^c$.
\end{thm}
This remains open, and has been proved only for a very limited set of graphs $H$ (although see~\cite{density4, density7} for a variety
of new graphs $H$ that satisfy \ref{EHconj}).

If $H$ is a graph, for each $n>0$ let $f_H(n)$ be the largest integer such that every $H$-free
 graph with at least $n$ vertices has a stable set or clique with size at least $f_H(n)$.
Thus, the \erh{} conjecture says that
\begin{thm}\label{EHconj2}
{\bf Conjecture: }For every graph $H$ there exists $c>0$ such that $f_H(n)\ge n^c$ for all $n\ge 0$.
\end{thm}
What can we actually prove about $f_H(n)$?
Erd\H{o}s and Hajnal~\cite{EH89} proved (logarithms are to base two, throughout the paper):
\begin{thm}\label{EHthm}
For every graph $H$, there exists $c>0$ such that $f_H(n)\ge 2^{c\sqrt{\log n}}$ for all $n> 0$.
\end{thm}

In~\cite{density1}, with Buci\'c, we improved this, showing:
\begin{thm}\label{loglog}
For every graph $H$, there exists $c>0$ such that $f_H(n)\ge 2^{c\sqrt{\log n\log \log n}}$ for all $n> 0$.
\end{thm}
For general graphs $H$, this is the best bound known.

The case when $H$ is a path is of particular interest. Paths are very simple graphs, and yet
until recently, the \erh{} conjecture was open when $H$ is the five-vertex path $P_5$. This was the smallest graph for which the conjecture was open,
and was therefore a focus of attention.
In a substantial breakthrough for the case $H=P_5$, Blanco and Buci\'c~\cite{blbu} improved the power in the exponent from $1/2$ to $2/3$: 
\begin{thm}\label{blbu}
There exists $c>0$ such that $f_{P_5}(n)\ge 2^{c(\log n)^{2/3}}$ for all $n> 0$.
\end{thm}
Their argument uses a complex structural analysis of $P_5$-free graphs, and does not appear to extend to longer paths.  However,
more recently, using completely different methods, we fully resolved the case $H=P_5$, showing in~\cite{density7} that $P_5$ satisfies the \erh{} conjecture:
\begin{thm}
There exists $c>0$ such that $f_{P_5}(n)\ge n^c$ for all $n>0$.
\end{thm}
But the \erh{} conjecture
for general paths remains open,
and it is natural to ask whether it is possible to obtain an improvement in the exponent such as \ref{blbu} for general paths.

In the present paper, we do better than that: 
we show that paths of any length have the ``near \erh{}'' property, that is, for every $d<1$ and every path $H$, there exists $c>0$ such that 
every $H$-free graph $G$ has  a clique or 
stable set of size at least $2^{(c\log |G|)^d}$ (note that $d=1$ is the \erh{} conjecture itself).  In other words:
\begin{thm}\label{mainthm}
For every path $P$, $f_{P}(n)\ge 2^{(\log n)^{1-o(1)}}$.
\end{thm}

We prove a stronger form of this (\ref{mainrodlthm} below), showing that $P$-free graphs contain large induced subgraphs which are very sparse or very dense.  Before stating this, we need some further discussion about dense and sparse subgraphs.

For $\vare>0$, a subset $S\subseteq V(G)$ is
\begin{itemize}
\item {\em $\vare$-sparse}
if the induced subgraph $G[S]$ has maximum degree at most $\vare|S|$;
\item {\em $(1-\vare)$-dense} if $\overline{G}[S]$ is $\vare$-sparse, where
$\overline{G}$ is the complement graph of $G$; and
\item {\em $\vare$-restricted} if $S$ is either $\vare$-sparse or $(1-\vare)$-dense.\footnote{Note that we are using  `local' density conditions: for $\epsilon$-sparse sets, we demand that every vertex has degree at most $\vare|S|$, rather than the `global' condition that there are at most $\vare\binom{|S|}{2}$ edges; and similarly for $(1-\epsilon)$-dense sets.  It is straightforward to move between the two, at the cost of a constant factor in the value of $\vare$ and the size of $S$.}
\end{itemize}
Note that if $|S|\le2$ then $S$ is trivially $\epsilon$-restricted.

An important result of R\"odl~\cite{rodl} shows that, for fixed $\epsilon>0$, $H$-free graphs have $\epsilon$-sparse or $(1-\epsilon)$-dense subsets of linear size:
\begin{thm}\label{rodlthm}
For all $0<\vare\le 1/2$, there exists $\delta>0$ such that for every $H$-free graph $G$,
there is an $\vare$-restricted subset $S\subseteq V(G)$ with $|S|\ge \delta|G|$.
\end{thm}
This does not imply the Erd\H os-Hajnal conjecture, as the dependence of $\delta$ on $\epsilon$ could be very poor. Indeed, R\"odl's proof used Szemer\'edi's Regularity Lemma, which leads to tower-type dependence. 
However, Fox and Sudakov~\cite{foxsudakov} conjectured that a {\em polynomial} version of R\"odl's theorem should hold, where $\delta$ depends polynomially on $\vare$.  More exactly:
\begin{thm}\label{foxconj}
{\bf Fox-Sudakov Conjecture: }For every graph $H$ there exists $c>0$ such that for every $\vare$ with $0<\vare\le 1/2$ and every $H$-free graph $G$,
there exists $S\subseteq V(G)$
with $|S|\ge \vare^c|G|$ such that $S$ is $\vare$-restricted.
\end{thm}
Let us say that a graph $H$ is {\em polynomial R\"odl} if it satisfies \ref{foxconj}.
It is straightforward to show that if $H$ is polynomial R\"odl then $H$ also satisfies the \erh{} conjecture.  
Indeed, the Fox-Sudakov conjecture interpolates between R\"odl's theorem and the \erh{} conjecture: taking $\vep$ to be a constant gives \ref{rodlthm}, and taking $\vep$ to be a small negative power of $n$ yields \ref{EHconj}.
Recently, 
Buci\'c, Fox and Pham \cite{bfp} showed the reverse implication:~the \erh{} conjecture and the Fox-Sudakov conjecture are in fact equivalent.

Both conjectures are currently out of reach.  But the bound on $\delta$ in \ref{rodlthm} has been strengthened:
Fox and Sudakov~\cite{foxsudakov} proved that \ref{rodlthm} holds with
$$\delta=\vep^{c_H\log(1/\vep)},$$
and this was improved in \cite{density1} to
$$\delta=\vep^{c_H\frac{\log(1/\vep)}{\log\log(1/\vep)}}.$$

In order to prove \ref{mainthm}, we will show that a much better bound than this holds when $H$ is a path: 
we will prove that every path ``nearly'' satisfies the Fox-Sudakov conjecture. 
Here is our main result:
\begin{thm}\label{mainrodlthm}
For every path $P$ and every $\alpha>0$, there exists $c>0$ such that for all $0<\vare\le 1/2$, every $P$-free graph $G$ contains an $\vare$-restricted subset $S\subseteq V(G)$ with $|S|\ge \delta|G|$, where
$$\delta=\vep^{c(\log \frac 1\vep)^{\alpha}}.$$
\end{thm}
Let us show that \ref{mainrodlthm} implies \ref{mainthm}.

\begin{thm}\label{restrictedtoclique}
Let $G$ be a graph with $|G|\ge 2$, let $a>0$ and $c:=(2a+2)^{-1}$. 
Let $\alpha\ge 0$, and assume that, for every $\vep\in(0,1/2)$,
there is an $\vep$-restricted $S\subseteq V(G)$ with
$$\abs S\ge \vep^{a(\log\frac1\vep)^{\alpha}}|G|.$$
Then $G$ contains a clique or stable set of size at least
$2^{c(\log |G|)^{\beta}}$, where $\beta=1/(1+\alpha)$.
\end{thm}
\Proof
Since $G$ has a clique or stable set of size two, 
we may assume that
$2^{c(\log |G|)^{\beta}}>2$.  Set
$$\vep:=2^{-2c(\log |G|)^{1/(1+\alpha)}}$$ 
and let
$$
\delta
=\vep^{a(\log\frac 1\vare)^{\alpha}}
=2^{-a(\log\frac 1\vare)^{1+\alpha}}
=2^{-a(2c)^{1+\alpha}\log|G|}
\ge|G|^{-2ac}.
$$
From the hypothesis, there is an $\vep$-restricted $S\subseteq V(G)$ with
 \[\abs S\ge\delta\abs G\ge |G|^{1-2ac}=|G|^{2c}
=2^{2c\log |G|}\ge 2^{2c(\log |G|)^{\beta}}
=\vep^{-1}.\]
Thus, since $S$ is $\vep$-restricted, $G[S]$ or its complenent has maximum degree at most $\vep|S|$ and so contains a clique or stable set of size at least
\[\frac{\abs S}{\vep\abs S+1}
\ge\frac1{2\vep}
\ge\vep^{-1/2}
= 2^{c(\log |G|)^{\beta}}.\]
This proves \ref{restrictedtoclique}.~\bbox

Note that if the Fox-Sudakov conjecture held, then we could take $\alpha=0$ in \ref{mainrodlthm}.  We could then choose $\epsilon$ to be some small negative power of $|G|$ (so that $\log(1/\vep)\approx\log |G|$), and applying \ref{mainrodlthm} would then give an $\vare$-restricted subset of polynomial size.  By Tur\'an's theorem or a greedy argument, this would yield a clique or stable set of polynomial size, giving the Erd\H os-Hajnal theorem for $P$-free graphs.  However we do not quite have this: since we must take $\alpha>0$ we do not quite get polynomial dependence, as if we take $\log(1/\vare)\approx \log n$ then the set $S$ given by \ref{mainrodlthm} is too small.  Instead, we need to take $\log(1/\vare)\approx(\log n)^{\beta}$ in the proof of \ref{restrictedtoclique}, where $\beta=1/(1+\alpha)$ is slightly less than 1. 

The rest of the paper is organized as follows.  In section \ref{blockades}, we discuss the relationship between dense or sparse multipartite structures (blockades) and dense or sparse subgraphs.  We discuss first the (realtively straightforward) sparse case and then the (harder) dense case in section \ref{sparsecase}, and we put the proof of \ref{mainrodlthm} together in section \ref{densityargument}.  We conclude with some further discussion in section \ref{conclusion}.

We use standard terminology throughout.  If $G$ is a graph, $G[X]$ denotes the
induced subgraph with vertex set $X$; $|G|$ denotes the number of vertices of $G$; and $\overline{G}$ is the complement graph of $G$.

\section{Blockades}\label{blockades}

Recall that a graph $H$ is {\em polynomial R\"odl} if there exists $c>0$ such that for all $\vare\in(0,1/2)$, every $H$-free 
graph $G$ contains an $\vare$-restricted subset $S\subseteq V(G)$ with $|S|\ge\delta|G|$, where
$$\delta=\vare^c|G|.$$
We will say that a graph $H$ is 
{\em near-polynomial R\"odl} if for every $\alpha>0$, there exists $c>0$ such that for every $\vep\in(0,\frac12)$, 
every $H$-free graph $G$ contains an $\vep$-restricted $S\subseteq V(G)$ with $\abs S\ge \delta\abs G$, where
$$\delta=\vep^{c(\log (1/\vare))^{\alpha}}.$$
Thus, our main theorem \ref{mainrodlthm} says that every path is near-polynomial R\"odl.  

In order to prove \ref{mainrodlthm}, we will need to find dense or sparse multipartite structures which we refer to as blockades. A {\em blockade} in a graph $G$ is
        a finite sequence $(B_1,\ldots,B_n)$ of (possibly empty) disjoint subsets of $V(G)$;
        its {\em length} is $n$ and its {\em width} is $\min_{i\in[n]}\abs{B_i}$.
        For $k,w\ge0$, $(B_1,\ldots,B_n)$ is a {\em $(k,w)$-blockade} if its length is at least $k$ and its width is at least $w$.

We will need certain density conditions to hold.
If $A,B\subseteq V(G)$ are disjoint, we say $A$ is {\em $x$-sparse} to $B$ if every vertex in $A$ has at most $x|B|$
neighbours in $B$, and $A$ is {\em $(1-x)$-dense } to $B$ if every vertex in $A$ has at least $(1-x)|B|$ neighbours in $B$.
        For $x\in(0,\frac12)$, a blockade $(B_1,\ldots,B_n)$ is {\em $x$-sparse} if $B_j$ is $x$-sparse to $B_i$ for all $i,j\in[n]$ with $i<j$,
        and {\em $(1-x)$-dense} if $B_j$ is $(1-x)$-dense to $B_i$ for all $i,j\in[n]$ with $i<j$ (note the asymmetry: vertices in sets
of the blockade are dense or sparse to earlier sets in the blockade, but not necessarily to later sets).  We say that a blockade is {\em complete} if, for all $i\ne j$, every vertex of $B_i$ is adjacent to every vertex of $B_j$; {\em anticomplete} if (for all $i\ne j$) every vertex of $B_i$ is nonadjacent to every vertex of $B_j$; and {\em $x$-restricted} if it is $x$-sparse or $(1-x)$-dense.

There are three parameters we care about: the length, width, and sparsity (or density).
It is easier to prove that certain graphs contain blockades with some desired combination of the three parameters, than to 
prove directly that they contain large $\vare$-restricted sets. But the reason blockades are useful for us
is that there is a transference theorem (\ref{thm:trans} below), that says that if a graph $G$
and all its large induced subgraphs admit blockades with certain parameters, then $G$ must contain a large 
$\vare$-restricted set.

A function $\ellh\colon(0,\frac12)\to\mab R^+$ is {\em subreciprocal} if it is nonincreasing and $1<\ellh(x)\le 1/x$ for all $x\in(0,\frac12)$.
        For a subreciprocal function $\ellh$,
        a graph $H$ is {\em $\ellh$-dividing}\footnote{In previous papers, we use the term {\em $\ellh$-divisive} to refer to graphs having this property when a few copies of $H$ are allowed. In this paper, we work with the slightly simpler property of being $\ellh$-dividing.  Buci\'c, Fox and Pham \cite{bfp} showed that these two properties are equivalent.} 
        if there are $c\in(0,\frac12)$ and $d>1$ such that for every $x\in(0,c)$ every $H$-free graph $G$ contains an $x$-sparse or $(1-x)$-dense $(\ellh(x),\floor{x^d\abs G})$-blockade in $G$.
Here is the transference theorem that we will use (see~\cite{density1}, or~\cite{density4} for an alternative proof of a weaker statement).\footnote{The version in \cite{density1} uses edge density rather than $\vare$-restricted subsets.  However, it is straightforward to adjust for this by changing the constants.}

\begin{thm}
        \label{thm:trans}
        Let $\ellh$ be subreciprocal, and let $H$ be an $\ellh$-dividing graph.
        Then there exists $C>0$ such that for every $\vep\in(0,\frac12)$, 
        every $H$-free graph $G$ contains an $\vep$-restricted $S\subseteq V(G)$ with $\abs S\ge \delta\abs G$,
        where
                $$\delta=\vep^{C\log(\frac1\vep)/\log\ellh(\vep)}.$$
\end{thm}

This result allows us to move between dense or sparse multipartite substructures (blockades) and dense or sparse subsets.  This approach has been used a number of times: the challenge is to find blockades that are both long and wide, and satisfy strong restriction properties.  For general graphs $H$, it is possible to find blockades that very restricted but not very long: 
\begin{itemize}
\item Erd\H{o}s and Hajnal~\cite{EH89} proved that for every graph $H$, there exists $d>0$
such that 
for all $x\in (0,1/2]$, every $H$-free graph admits disjoint sets $A,B$ of size at least $\lfloor x^d|G|\rfloor$  such that $B$ is $x$-sparse or $(1-x)$-dense to $A$; in other words, 
an $x$-sparse or $(1-x)$-dense
$(2, \lfloor x^d|G|\rfloor)$ blockade. From this they deduced their result \ref{EHthm}.
\item In \cite{density1}, we (with Buci\'c) proved a strengthening: for every graph $H$, there exists $d>0$
such that
for all $x\in (0,1/2]$, every $H$-free graph
admits an $x$-sparse or $(1-x)$-dense
$(\log(1/x), \lfloor x^d|G|\rfloor)$-blockade.  This allowed us to deduce \ref{loglog}.
\end{itemize}

To prove results with {\em polynomial} dependence, blockades with much better parameters are required.  In general, it is a very hard challenge to find suitable blockades, but there has been some progress.  Let $\mathcal G$ be a class of graphs that is closed under taking induced subgraphs.

\begin{itemize}
\item If we can find complete or anticomplete pairs of {\em linear} size, then we can obtain polynomial-size cliques or stable sets.  More precisely, we say that $\mathcal G$ has the {\em strong Erd\H os-Hajnal property} if there is $c>0$ such that every $G\in \mathcal G$ admits a complete or anticomplete $(2,\lfloor c|G|\rfloor)$-blockade.  It is not hard to show that this is enough to give cliques or stable sets of polynomial size.  Unfortunately, for almost all $H$, the class of $H$-free graphs does not have the strong Erd\H os-Hajnal property, but there are interesting classes where it does hold (for example if we exclude both a forest and the complement of a forest \cite{pure1}).
\item It may be possible to trade off between length and width.  Suppose that there is some $d>0$ such that every $G\in\mathcal G$ with at least two vertices admits a $(k,|G|/k^d)$-blockade for some $k\ge 2$ (which may be different for different $G$) .  This is the ``quasi-Erd\H os-Hajnal property'' \cite{pt,tomon}, and is again sufficient to give cliques or stable sets of polynomial size.  It was shown in \cite{c5} that the class of $C_5$-free graphs has the quasi-Erd\H os-Hajnal property, proving the Erd\H os-Hajnal conjecture for $C_5$.
\item We can sometimes also use the density of $G$.  Suppose that every $G\in\mathcal G$ contains either a suitable blockade or a large induced subgraph with much better density parameters: say, every $y$-restricted $G\in\mathcal G$  contains either a suitably restricted $({\rm poly}(1/y), \lfloor {\rm poly}(y)|G|\rfloor)$-blockade, or a $\lambda y$-restricted subset of size at least ${\rm poly}(\lambda)|G|$, where $\lambda\in(0,1)$ may depend on $G$.  Iterating this, we can pass to a succession of ever more restricted subgraphs until we reach either a polynomial size clique or stable set or a suitable blockade.  This is the method of iterative sparsification (introduced in \cite{density4}, and further developed in \cite{density6,density7}).
\end{itemize}
In all these results, the main obstacle is proving the existence of blockades that are sufficiently long, wide and restricted.

In this paper we are interested in $P$-free graphs, where $P$ is a path.  If we could show that $P$ is $\ellh^*$-dividing, where 
$$\ellh^*(x)=(1/x)^c=2^{c\log\frac1x},$$
then \ref{thm:trans} would yield the polynomial R\"odl property.  However, the blockades that we obtain are not sufficiently long to show this.  Instead, we will will work with a sequence of non-polynomial approximations to $h^*$.  
For each integer $s\ge 0$, let 
$\ellh_s\colon(0,\frac12)\to\mab R^+$ be the function defined by
        $$\ellh_s(x):=2^{(\log\frac1x)^{\frac{s}{s+1}}}$$
        for all $x\in(0,\frac12)$.
        (Thus $\ellh_s$ is subreciprocal.)
We will show that:
\begin{thm}
        \label{thm:main}
Every path $P$ is  $\ellh_s$-dividing for all integers $s\ge0$.
\end{thm}

Let us first deduce our main result \ref{mainrodlthm}.  

\medskip

\noindent
{\bf Proof of \ref{mainrodlthm} (assuming \ref{thm:main}).\ \ }
We must show that for 
all $\alpha>0$, there exists $c>0$ such that for all $\vare\in (0, 1/2)$, every $H$-free graph $G$ contains an $\vep$-restricted $S\subseteq V(G)$ with $\abs S\ge \delta\abs G$, where
$$\delta=\epsilon^{c(\log\frac{1}{\vep})^\alpha}.$$

Choose an integer $s$ such that $\frac{1}{s+1}\le \alpha$. Since $P$ is $\ellh_s$-dividing, by \ref{thm:trans} there exists $C>0$ such that 
for every $\vep\in(0,\frac12)$, every $H$-free graph $G$ contains an $\vep$-restricted $S\subseteq V(G)$ with $\abs S\ge \eta\abs G$, where
$$
\eta
=\vep^{C\log(\frac{1}{\vep})/\log(\ellh_s(\vep))}
=\vep^{C\log(\frac{1}{\vep})/(\log\frac{1}{\vep})^{\frac{s}{s+1}}}
=\vep^{C(\log(\frac{1}{\vep}))^{\frac{1}{s+1}}}
\ge\vep^{C(\log(\frac{1}{\vep}))^{\alpha}}.
$$
Taking $c=C$, we have $\eta\ge\delta$, and the result follows. This proves \ref{mainrodlthm}.~\bbox
As noted in the previous section, \ref{mainrodlthm} implies that every path is near-polynomial R\"odl.  

The proof of \ref{thm:main} will run by induction on $s$: to show that $P$ is $\ellh_s$-dividing, we first use the fact that $P$ is $h_{s-1}$-dividing to find a large restricted subset, and then work inside this to obtain the required blockade.
We will need to deal with two cases: an easier sparse case and a harder dense case.  We address these separately in the next section.  The main argument is then given in section \ref{densityargument}.

\section{Sparse and dense graphs}\label{sparsecase}

The proof of \ref{thm:main} proceeds by induction on $s$.  Every graph is $\ellh_0$-dividing, 
so we may assume that $s>0$ and that the path $P$ is $\ellh_{s-1}$-dividing. 
We are given a $P$-free graph $G$, and to show
that $P$ is $\ellh_{s}$-dividing we need to find a blockade in $G$ with certain parameters. By \ref{thm:trans}, we know that every large subset of $V(G)$ includes a somewhat smaller subset $S$ that is either very  dense or very 
sparse.

There is no symmetry between the dense and sparse cases, because we are excluding a path but not its complement.  As usual with 
problems about excluding a path, our task is easier if the host graph is sparse (see, for example, \cite{pure1}).  
We first prove a result that will allow us to win if ever the subset $S$ provided by \ref{thm:trans} turns out to be sparse.  

\begin{thm}\label{newsparse}
Let $P$ be a path with $k\ge2$ vertices, and let $0<y\le 1/(60k)$.  Let $G$ be a $y^2$-sparse graph. Then either:
\begin{itemize}
\item $G$ contains a copy of $P$, or
\item $G$ contains an anticomplete $(1/y,\lfloor y^{2}|G|\rfloor)$-blockade in $G$.
\end{itemize}
\end{thm}

\Proof
We may assume that $|G|\ge 1/y^2$ or the result is trivial. Let $\mu=1/(12y)$.  

Let us say that a set $A$ of vertices in a graph $G$ is {\em poorly expanding} (in $G$) if $|N_G(A)|\le\mu|A|$, where $N_G(A)$ is the set of vertices $u\not\in A$ that have a neighbour in $A$. We repeatedly remove poorly expanding sets together with their neighbourhoods, as follows: let $G_0=G$ and, for $i\ge0$, if $G_i$ contains a set $A_{i+1}$ of vertices of size at most $y^2|G|$ that is poorly expanding in $G_i$, then let $B_{i+1}=N_{G_i}(A_{i+1})$ and set $G_{i+1}=G_i\setminus(A_{i+1}\cup B_{i+1})$.  The process terminates with some induced subgraph $G_t$ that has no poorly expanding set of size at most $y^2|G|$.  

Suppose first that $|G_t|\le|G|/2$.
Let $A=\bigcup_{i=1}^t A_i$ and $B=\bigcup_{i=1}^t B_i$, so $G_t=G\setminus(A\cup B)$ and $|B|\le\mu|A|$.  Since $|A\cup B|\ge|G|/2$, it follows that $|A|\ge|G|/(2+2\mu)>3y|G|$.  The sets $A_i$ are pairwise anticomplete, so each component of $G[A]$ has size at most $y^2|G|$.  We group the components together into as many sets as possible such that each contains between $y^2|G|$ and $2y^2|G|$ vertices.  If we get at least $1/y$ sets then we have found our blockade.  If not, then there are at most $y^2|G|$ unused vertices and so $|A|\le(1/y)2y^2|G|+y^2|G|<3y|G|$, giving a contradiction.

Otherwise, $|G_t|>|G|/2$ and $G_t$ contains no poorly expanding set of size at most $y^2|G|$.  We work inside $G_t$ for the rest of 
the argument.  We define sets of vertices $R_1,\dots,R_{k-1}$ and $S_0,\dots,S_{k-1}$ as follows: let $S_0\subseteq V(G_t)$ 
be any set of $\lfloor y^2|G|\rfloor$ vertices of $G_t$; for $i\ge1$, let $R_i$ be a minimal subset of $S_{i-1}$ such that 
$|N_{G_t}(R_i)\setminus \bigcup_{j<i}S_j|\ge y^2|G|$; and let $S_i=N_{G_t}(R_i)\setminus \bigcup_{j<i}S_j$.  If we can do this 
successfully, then we can find an induced copy $v_1\CC v_k$ of $P$: let $v_k$ be any vertex of $S_{k-1}$ and, for $i=k-1,\dots,1$, let $v_i$ be any neighbour of $v_{i+1}$ in $R_i$.

It is therefore enough to show that we can find sets $R_i$ and $S_i$ at each stage.  Suppose that $1\le i\le k-1$ and we have defined $S_0,\dots,S_{i-1}$.  
Let $T\subseteq S_{i-1}$ be any subset of size $\lfloor y^2|G|\rfloor$.  Then, by our assumption on poorly expanding sets, 
$|N_{G_t}(T)|\ge \mu\lfloor y^2|G|\rfloor \ge y|G|/24$.  By the minimality of $R_j$ and since $G$ is $y^2$-sparse, 
each set $S_j$ has size at most $2y^2|G|$ and so $|N(T)\setminus\bigcup_{0\le j< i}S_j|\ge y|G|/24 - 2ky^2|G|\ge y^2|G|$.  
It follows that we can pick $R_{i}$ and $S_{i}$.
This proves \ref{newsparse}.~\bbox

The dense case is much more difficult.
This is handled by the next lemma, which will allow us either to extract a dense blockade, or to move to a reasonably large induced subgraph with no large, dense parts.

\begin{thm}\label{newdense}
Let $P$ be a path with $|P|\ge 1$, let $a=3|P|$, and let $0<x\le y\le 1/100$. For every $P$-free graph $G$,
either:
\begin{itemize}
\item[(a)] there is $S\subseteq V(G)$ 
with $|S|\ge x^{a}|G|$ such that every $S'\subseteq S$ with $|S'|\ge x|S|$ has density at most $(1-y^3)$; or
\item[(b)] there is a $(1-x)$-dense $(1/y,\lfloor x^{a}|G|\rfloor)$-blockade in $G$.
\end{itemize}
\end{thm}
\Proof
Let us assume that neither (a) nor (b) holds.
It follows that $x^{a}|G|\ge 1$, since otherwise (b) holds.
A {\em non-neighbour} of a vertex $v$ means a vertex different from and nonadjacent to $v$, and the
{\em antidegree} of $v$ is the number of its non-neighbours.
We claim first:
\\
\\
(1) {\em For every $S\subseteq V(G)$ with $|S|\ge 2x^{a-1}|G|$, there exists $C\subseteq S$ with $|C|\ge x(1+y)|S|/2$, such that
\begin{itemize}
    \item $C$ is $(1-2y^3)$-dense; and
    \item  for all disjoint $X,Y\subseteq C$ with $|X|\ge (1-y/4)|C|$ and $|Y|\ge x^a|G|$, 
at least $y|X|/4$ vertices in $X$ have at least $x|Y|$ non-neighbours in $Y$.
\end{itemize}
 }
\noindent Since $|S|\ge x^a|G|$, and (a) does not hold, there exists $S'\subseteq S$ with $|S'|\ge x|S|$ such that $S'$ is $(1-y^3)$-dense. 
Choose a $(1-x)$-dense blockade $(B_1\LL B_{n-1},C)$ in $G[S']$ with $n$ maximum such that $|B_1|\LL |B_{n-1}|,|C|\ge x^a|G|$ and 
$|C|\ge (1-(n-1)y/2)|S'|$. (This is possible because $|S'|\ge x^a|G|$, and so we can take $n=1$ and $C=S'$.) We may assume that $n<1/y$, and so 
$$|C|\ge (1-(n-1)y/2)|S'|= (1+y/2-ny/2)|S'|\ge (1+y) |S'|/2\ge x(1+y)|S|/2.$$
In particular, $|C|\ge |S'|/2$, and consequently $C$ is $(1-2y^3)$-dense.
\begin{figure}[h!]
\centering

\begin{tikzpicture}[scale=1,auto=left]

\tikzstyle{every node}=[]
\draw (0,0) ellipse (6 and 1);
\draw (0,0) ellipse (5 and .9);
\draw[rounded corners]  (-4.2,-.4) rectangle (-3.4,.4);
\draw[rounded corners]  (-2.3,-.4) rectangle (-1.5,.4);
\draw[rounded corners]  (-1.2,-.5) rectangle (4,.5);
\draw[dotted, very thick] (-3.3,0)--(-2.4,0);

\node at (-5.5,0) {$S$};
\node at (-4.5,0) {$S'$};
\node at (-3.8,0) {$B_1$};
\node at (-1.9,0) {$B_{n-1}$};
\node at (1.4,0) {$C$};


\end{tikzpicture}

\caption{For step (1) of the proof of \ref{newdense}.} \label{fig:useblockade}
\end{figure}

Suppose that $X,Y\subseteq C$ are disjoint, with $|X|\ge (1-y/4)|C|$ and $|Y|\ge x^a|G|$. It follows that 
$$|X|\ge (1-y/4)(1+y)|S'|/2\ge |S'|/2.$$
Since $|Y|\ge x^a|G|$,
and $(1-ny/2)|S'|\ge |S'|/2\ge x^a|G|$, 
fewer than $(1-ny/2)|S'|$ vertices in $X$ are $(1-x)$-dense to $Y$, from the maximality of $n$. Since $|C|\ge (1-(n-1)y/2)|S'|$,
it follows that at least $y|S'|/2-|Y|$ vertices in $X$ have at least $x|Y|$ non-neighbours in $Y$.
But $|Y|\le y|C|/4$, since $X\cap Y=\emptyset$, and so 
$y|S'|/2-|Y|\ge y|S'|/2- y|C|/4\ge y|C|/4\ge y|X|/4$. This proves (1).

\bigskip
Next we use a new version of the ``Gy\'arf\'as path'' argument, originally used to find induced paths in graphs with large chromatic
number (see \cite{gy1,gy2}).  The usual form of the argument grows a path $v_1\CC v_k$ one vertex at a time, walking towards some
part of the graph with large chromatic number: at each stage, $v_i$ has neighbours in some connected part of the graph with large
chromatic number that is nonadjacent to $\{v_1,\dots,v_{i-1}\}$, and we choose $v_{i+1}$ from this part.  If the neighbourhood of
every vertex has small chromatic number, this can be made to work.  However, the current setup is rather different: we are dealing
with cardinality rather than chromatic number, and we are working inside a graph that is very dense.  There is no problem arranging
that the last vertex
of the path has many neighbours; the issue is to arrange that there are many vertices with no neighbours in the path, and to maintain this as we grow the path.

Let $k:=|P|$, so $a=3k$. 
Define $a_1:=x/2$, and $b_1:=x^2y/8$; and for $2\le t\le k$, define $a_t := x^{2t-1}y/2^{t+2}= (x/2)b_{t-1}$ and $b_t:=x^{2t}y/2^{t+2}=(x^2/2)b_{t-1}$.
For $1\le t\le k$ let us say a {\em $t$-brush} is an induced path of $G$ with vertices $v_1\CC v_t$
in order, such that there exist subsets $A,B\subseteq V(G)\setminus \{v_1\LL v_t\}$ with the following properties:
\begin{itemize}
\item every vertex in $A$ is adjacent to $v_t$ and is nonadjacent to $v_1\LL v_{t-1}$;
\item every vertex in $B$ has no neighbours in $\{v_1\LL v_t\}$;
\item $|A|\ge a_t|G|$ and $|B|\ge b_t|G|$;
\item for every $Y\subseteq B$ with $|Y|\ge x^a|G|$, there are at least $y|A|/4$ vertices in $A$ that have at least $x|Y|$
non-neighbours in $Y$;
and
\item every vertex in $B$ has at most $3y^3|A|$ non-neighbours in $A$.
\end{itemize}
\noindent 
(2) {\em There is a $1$-brush.}
\\
\\
Since $|G|\ge 2x^{a-1}|G|$, (1) implies that there exists $C\subseteq V(G)$ with $|C|\ge x(1+y)|G|/2$, such that $C$ is 
$(1-2y^3)$-dense, and
and for all disjoint $X,Y\subseteq C$ with $|X|\ge (1-y/4)|C|$ and $|Y|\ge x^a|G|$,
at least $y|X|/4$ vertices in $X$ have at least $x|Y|$ non-neighbours in $Y$.

Suppose first that no vertex of $C$ has at least $b_1|G|$ nonneighbours in $C$.  
Let $Y\subseteq C$ be a set
of $\lceil x^a |G|\rceil$ vertices and let $X=C\setminus Y$.  
Then the total number of nonedges between $X$ and $Y$ is at most the sum of antidegrees of vertices in $Y$, which is at most $b_1|G| |Y|/2=x^2y|G| |Y|/8$.  On the other hand, 
$$|Y|\le x^a |G| + 1\le 2x^a |G|\le 4x^{a-1}|C|\le (x/4)|C|\le (y/4)|C|,$$
so by our choice of $C$ the number of nonedges is at least $(y|X|/4)x|Y|$.  It follows that
$$(y|X|/4)x|Y|\le x^2y|G| |Y|/8,$$ 
and so $|X|\le x|G|/2$, which gives a contradiction as $|X|\ge(1-y/4)|C|> x|G|/2$.

So let $v_1$ be a vertex of $C$ with at least $b_1|G|$ nonneighbours in $C$.  Let $A$ be the set of its nonneighbours in $C$, and 
let $B=C\setminus A$.  As $C$ is $(1-2y^3)$-dense, we have 
$$|A|\ge (1-2y^3)|C|\ge x|G|/2=a_1|G|;$$
and we clearly have $|B|\ge b_1|G|$.  Since $A$ and $B$ are subsets of $C$, the last two bullets in the definition of a $1$-brush are also satisfied. This proves (2).
\\
\\
(3) {\em Let $1\le t\le k-1$, and let $v_1\CC v_t$ be a $t$-brush. Then there is a vertex $v$ such that
$v_1\CC v_t\DD v$ is a $(t+1)$-brush.}
\\
\\
Choose $A,B$ satisfying the five bullets in the definition of ``$t$-brush''. 
Since $b_t=(x^2/2)^ty/4$, and $t\le k-1$, and $a\ge 3k$, it follows that 
$$|B|\ge b_t|G|= (x^2/2)^ty|G|/4\ge x^{3k-4}|G| \ge 2x^{a-1}|G|.$$
By (1), there exists $C\subseteq B$ with $|C|\ge x(1+y)|B|/2$, such that
$C$ is $(1-2y^3)$-dense, and for all disjoint $X,Y\subseteq C$ with $|X|\ge (1-y/4)|C|$ and $|Y|\ge x^a|G|$,
at least $y|X|/4$ vertices in $X$ have at least $x|Y|$ non-neighbours in $Y$.

Since $v_1\CC v_t$ is a $t$-brush, each vertex in $C$ has at most $3y^3|A|$ non-neighbours in $A$, and so at most $y|A|/8$
vertices in $A$ have at least $24y^2|C|$ non-neighbours in $C$. On the other hand, 
there are at least $y|A|/4$ vertices in $A$ that have at least $x|C|$
non-neighbours in $C$; and so there is a set $D\subseteq A$ with $|D|\ge y|A|/8$, such that for each $v\in D$,
the number of its non-neighbours in $C$ is between $x|C|$ and $24y^2|C|$.

\begin{figure}[H]
\centering

\begin{tikzpicture}[scale=0.8,auto=left]
\tikzstyle{every node}=[inner sep=1.5pt, fill=black,circle,draw]
\node (v1) at (0,2.5) {};
\node (v2) at (0,1.5) {};
\node (vt) at (0,0) {};
\node (v) at (-2,-1.3) {};

\foreach \from/\to in {v1/v2}
\draw [-] (\from) -- (\to);

\draw[dotted, very thick] (0,.2) to (0,1.3);

\tikzstyle{every node}=[]
\draw (0,-1.2) ellipse (6 and .4);
\draw (0,-2.2) ellipse (4 and .4);
\draw[rounded corners]  (-2,-2.5) rectangle (2, -1.9);
\draw (-2.3,-1.2) ellipse (2 and .3);
\draw[rounded corners]  (-1.9,-2.4) rectangle (.9, -2);
\draw (vt) to (-5,-1.1);
\draw (vt) to (-3,-1.1);
\draw (vt) to (-1,-1.1);
\draw (vt) to (1,-1.1);
\draw (vt) to (3,-1.1);
\draw (vt) to (5,-1.1);

\draw (v) to (-1.8,-2.2);
\draw (v) to (-1.15,-2.2);
\draw (v) to (.15,-2.2);
\draw (v) to (.8,-2.2);

\draw[right] (v1) node []           {$v_1$};
\draw[right] (v2) node []           {$v_2$};
\draw[above right] (vt) node []           {$v_t$};
\node at (2.5,-1.2) {$A$};
\node at (2.5,-2.2) {$B$};
\node at (-3.5,-1.2) {$D$};
\node at (-2.3,-1.2) {$v$};
\node at (-.8,-2.2) {$A'$};
\node at (1.5,-2.2) {$B'$};

\end{tikzpicture}
\caption{For step (3). $C=A'\cup B'$.} \label{fig:step3}
\end{figure}

Let $v\in D$. We claim that $v_1\CC v_t\DD v$
is a $(t+1)$-brush. 
Let $A'$ be the set of all neighbours of $v$ in $C$, and let $B'=C\setminus A'$. 
We will show that $A',B'$ satisfy the five conditions in the definition of a $(t+1)$-brush.
The first two are immediate. For the third, 
$$|A'|\ge (1-24y^2)|C|\ge (1-24y^2)x(1+y)|B|/2\ge (x/2)b_t|G|= a_{t+1}|G|,$$
and 
$$|B'|\ge x|S|\ge x(x/2)|B|\ge (x^2/2)b_t|G|= b_{t+1}|G|.$$
For the fourth condition, suppose that 
$Y\subseteq B'$ with $|Y|\ge x^a|G|$. From the choice of $C$, since $|A'|\ge (1-24y^2)|C|\ge (1-y/4)|C|$, there are at 
least $y|A'|/4$ vertices 
in $A'$ that have at least $x|Y|$
non-neighbours in $Y$. Finally, for the fifth condition, since $C$ is $(1-2y^3)$-dense, each vertex in $B'$ has at most 
$2y^3|C|\le 3y^3|A'|$ non-neighbours in $A'$. This proves (3).

\bigskip

The result now follows by induction: (2) shows that there is a 1-brush, and (3) provides the inductive step.  Thus there is a $k$-brush, as required.  This proves \ref{newdense}.~\bbox

\section{Decreasing density}\label{densityargument}
We remind the reader that if $\ellh$ is subreciprocal, a graph $H$ is {\em $\ellh$-dividing} if there are 
$c\in(0,\frac12)$ and $d>1$ such that for every $x\in(0,c)$ and every $H$-free graph $G$ admits an $x$-sparse or $(1-x)$-dense $(\ellh(x),\floor{x^d\abs G})$-blockade in $G$. 
Also, for each integer $s\ge 0$, $\ellh_s\colon(0,\frac12)\to\mab R^+$ is the function defined by 
        $$\ellh_s(x):=2^{(\log\frac1x)^{\frac{s}{s+1}}}$$
        for all $x\in(0,\frac12)$. In this section we use the results of the previous two sections, together with 
\ref{thm:trans}, to prove \ref{thm:main}, which we restate:

\begin{thm}
        \label{thm:main2}
Every path $P$ is  $\ellh_s$-dividing for all integers $s\ge0$.
\end{thm}
\Proof
The proof is by induction on $s$.  The key part of the argument is the inductive step, where we bootstrap from $h_{s-1}$-dividing to $h_s$-dividing.  Our goal is to find a large blockade that is suitably dense or sparse.  As $P$ is $h_{s-1}$-dividing, we can use the transference lemma \ref{thm:trans} to obtain a large subset $S$ that induces subgraph that is either dense or sparse.  If the induced subgraph is sparse, then we can use \ref{newsparse} to obtain the required blockade; otherwise, $S$ induces a large dense subgraph.  In fact, we can apply this argument to any large induced subgraph: thus either we win, or we find that every large induced subgraph of $G$ contains a large, dense part.  This will allow us to use \ref{newdense} to obtain a dense blockade.

So let $G$ be a $P$-free graph.
For $s=0$,
we need to find disjoint sets $B_1$ and $B_2$ of size at least $\floor{x^d\abs G}$ such that $B_2$ is 
$x$-sparse or $(1-x)$-dense to $B_1$.  This follows from a strengthened version of a result of Erd\H{o}s and Hajnal~\cite{EH89} given by Fox and Sudakov~\cite{foxsudakov} (see their Lemma 2.1). So, we may assume that $s\ge 1$, and $P$ is $\ellh_{s-1}$-dividing. 

By \ref{thm:trans},
with $\ellh=\ellh_{s-1}$, we deduce that
there exists $C>0$ such that for every $\vep\in(0,\frac12)$, every $P$-free graph $G'$
contains an $\vep$-restricted $S\subseteq V(G')$ with $\abs S\ge \delta\abs G'$, where
$$\delta
=\vep^{\frac{C\log\frac1\vep}{\log(\ellh_{s-1}(\vep))}}
=\vep^{C(\log\frac1\vep)^{1/s}}.$$
Let $b=3^{1+1/s}C$.
We deduce:
\\
\\
(1) {\em Let $0< x\le 1/2$ and let $y:=1/\ellh_s(x)$. Then every $P$-free graph $G'$ contains
a $y^3$-restricted subset $S\subseteq V(G')$ with $\abs S\ge x^{b}|G'|$.}
\\
\\
Since $x\le 1/2$, it follows that $\ellh_s(x)\ge 2$ and so $y^3\le y\le 1/2$.
Moreover,
$$\log\frac 1y=\log(\ellh_s(x))=\left(\log\frac 1x\right)^{s/(s+1)}.$$
Setting $\vare=y^3$, we deduce that every $P$-free graph $G'$
contains a $y^3$-restricted $S\subseteq V(G')$ with $\abs S\ge \delta|G'|$, where
$$\delta=\vep^{C(\log\frac1\vep)^{1/s}}=2^{-C(\log\frac{1}{\vep})^{1+1/s}}=2^{-3^{1+1/s}C(\log\frac{1}{y})^{1+1/s}}=2^{-b(\log\frac{1}{y})^{1+1/s}}=x^{b}.$$
This proves (1).

\bigskip

Now, let $d=3b|P|+b+5$, and choose $c>0$ with $c\le 1/2$, and sufficiently small that 
$$c^{b}\le \frac{1}{\ellh_{s}(c)}\le \min\left(\frac{1}{60|P|},\frac{1}{100}\right).$$
Let $x\in(0,c)$ and let $G$ be a $P$-free graph. 
We will show that there is an $x$-sparse or $x$-dense $(\ellh_s(x),\floor{x^d\abs G})$-blockade in $G$, and therefore that $P$ is $\ellh_s$-dividing. Suppose (for a contradiction) that there is no such blockade.
As before, let $y:=1/\ellh_s(x)$. We note that $y\le x\le c\le 1/100$.
\\
\\
(2) {\em For every $S\subseteq V(G)$ with $|S|\ge x^{d-b-5}|G|$, there
exists a $(1-y^3)$-dense subset
$S'\subseteq S$ with $|S'|\ge x^{b}|S|$.}
\\
\\
Suppose not. By (1) applied to $G[S]$, 
there is an $y^3$-sparse subset $S'\subseteq S$ with $|S'|\ge x^{b}\abs S$.
By \ref{newsparse} applied to $G[S']$, there is an $x$-sparse $(1/y,\lfloor x^{5}|S'|\rfloor)$-blockade in $G[S']$.
But then $G$ admits an anticomplete $(1/y,\lfloor x^{b+2}y^2|G|\rfloor)$-blockade: this is 
an $x$-sparse $(1/y,\lfloor x^{d}|G|\rfloor)$-blockade,
giving a contradiction.
This proves (2).

\bigskip

In particular, (1) implies that for every $S\subseteq V(G)$
with $|S|\ge x^{3b|P|}|G|$, there is a $(1-y^3)$-dense subset $S'\subseteq S$ with $|S'|\ge x^{b}|S|$,
since $x^{3b|P|}= x^{d-b-5}$.
But now, by \ref{newdense} 
with $x$ replaced by $x^{b}$, we deduce
that there is a $(1-x^{b})$-dense, and hence $(1-x)$-dense, $(1/y,\lfloor x^{3b|P|}|G|\rfloor)$-blockade in $G$.
As $d\ge 3b|P|$, this proves \ref{thm:main2}.~\bbox

\section{Concluding remarks}\label{conclusion}
We have proved that all paths are near-polynomial R\"odl; can we do the same for any other graphs? A {\em caterpillar} is a tree in which there
is a path that contains at least one end of every edge of the tree; and the proof of this paper for paths can be adapted without 
difficulty to show that all caterpillars are near-polynomial R\"odl. This leads naturally to asking whether all trees are near-polynomial R\"odl, but we 
have not been able to prove that; we cannot handle the tree obtained from 
the 3-claw $K_{1,3}$ by subdividing each edge once (this is the unique minimal tree that is not a caterpillar).
But there are four further graphs (two complementary pairs) 
that we can prove are near-polynomial R\"odl, all with six vertices and shown in the figure \ref{fig:6vertex}. The method of proof is similar. 
We omit all these proofs, which will appear in Tung Nguyen's thesis~\cite{tungthesis}. 

\begin{figure}[h]
                \centering

\tikzstyle{every node}=[inner sep=1.5pt, fill=black,circle,draw]
                \begin{tikzpicture}[scale=0.8,auto=left]
                        \node (a1) at (-2,0) {};
                        \node (a2) at (-.5,0) {};
                        \node (a3) at (1,0) {};
                        \node (a4) at (-2,1.5) {};
                        \node (a5) at (-0.5,1.5) {};
                        \node (a6) at (1,1.5) {};

                        \node (b1) at (2,0) {};
                        \node (b2) at (4,0) {};
                        \node (b3) at (6,0) {};
                        \node (b4) at (3,1.5) {};
                        \node (b5) at (5,1.5) {};
                        \node (b6) at (7,1.5) {};

                        \node (c1) at (8.4,.5) {};
                        \node (c2) at (9,1.5) {};
                        \node (c3) at (9.6,.5) {};
                        \node (c4) at (7.6,0) {};
                        \node (c5) at (9,2.4) {};
                        \node (c6) at (10.4,0) {};

                        \node (d1) at (11.4,0) {};
                        \node (d2) at (13,0) {};
                        \node (d3) at (14.6,0) {};
                        \node (d4) at (12.2,1.3) {};
                        \node (d5) at (13.8,1.3) {};
                        \node (d6) at (13,2.6) {};

                        \foreach \from/\to in {a1/a2, a2/a3,a4/a5,a5/a6,a2/a5,a3/a6,b1/b2,b2/b3,b4/b5,b5/b6,b1/b4,b4/b2,b2/b5,b5/b3,b3/b6, c1/c2,c1/c3,c2/c3,c1/c4,c2/c5,c3/c6, d1/d2,d2/d3,d1/d4,d2/d4,d2/d5,d3/d5,d4/d5,d4/d6,d5/d6}
\draw [-] (\from) -- (\to);
\end{tikzpicture}
                \caption{
                        Four near-polynomial R\"odl six-vertex graphs.
                        } \label{fig:6vertex}
        \end{figure}

Another way to strengthen the current result is by looking at ordered graphs. An {\em ordered graph} $G$ is a pair  $(G^\natural, \le_G)$, where
$G^\natural$ is a graph and $\le_G$ is a linear order of its vertex set. Induced subgraph containment for ordered graphs
is defined in the natural way, respecting the orders of both graphs. A {\em zigzag path} is an ordered graph $(G^\natural, \le_G)$
where $G^\natural$ is a path and the ordering is as in figure \ref{fig:zigzag}. The proof of this paper works for ordered graphs,
with minor adjustments,
showing that
every zigzag path is near-polynomial R\"odl (defining ``near-polynomial R\"odl'' for ordered graphs in the natural way). We omit the details.

\begin{figure}[H]
\centering

\begin{tikzpicture}[scale=0.8,auto=left]
\tikzstyle{every node}=[inner sep=1.5pt, fill=black,circle,draw]

\node (v1) at (0,0) {};
\node (v2) at (2,0) {};
\node (v3) at (4,0) {};
\node (v4) at (6,0) {};
\node (v5) at (8,0) {};
\node (v6) at (10,0) {};
\node (v7) at (12,0) {};
\node (v8) at (14,0) {};

\draw [-] (v1) to [bend left=20] (v8);
\draw [-] (v2) to [bend left=20] (v8);
\draw [-] (v2) to [bend left=20] (v7);
\draw [-] (v3) to [bend left=20] (v7);
\draw [-] (v3) to [bend left=20] (v6);
\draw [-] (v4) to [bend left=20] (v6);
\draw [-] (v4) to [bend left=20] (v5);

\end{tikzpicture}
\caption{A zigzag path} \label{fig:zigzag}
\end{figure}

\bigskip
\noindent{\bf Note.}  For the purpose of Open Access, the authors have applied a CC BY public copyright licence to any Author Accepted Manuscript (AAM) version arising from this submission.

\end{document}